\documentclass[12pt]{amsart}
\usepackage{amscd,amssymb}

\newtheorem{thm}{Theorem}[section]
\newtheorem{corol}[thm]{Corollary} 

\newtheorem{prop}[thm]{Proposition}
\theoremstyle{definition}

\theoremstyle{remark}

\numberwithin{equation}{section}

\renewcommand{\Bbb}{\mathbb}

\def\Q {{\Bbb Q}}

\def\C {{\Bbb C}}

\def\T{{\mathbf T}}
\def\Z {{\mathbb Z}}

\def\QED{\nobreak\quad\ifmmode\roman{Q.E.D.}\else{\rm Q.E.D.}\fi}

\def\a{\alpha}

\def\Si{\Sigma}

\def\o{\omega}

\def\sR{{\mathcal R}}

\def\sbs{\subset}

\newcommand{\liminv}{\varprojlim}

\begin{document}

\title
{A compact group which is not Valdivia compact}

\author[W. Kubi\'s]{Wies\l aw Kubi\'s}

\address{Department of Mathematics, York University, 4700 Keele Street,
Toronto, Ontario
M3J 1P3, Canada}

\curraddr{The Fields Institute for Research in Mathematical Sciences,
222 College Street, Toronto, Ontario M5T 3J1, Canada}

\email{kubis@cs.bgu.ac.il}

\author[V. Uspenskij]{Vladimir Uspenskij}

\address{Department of Mathematics, 321 Morton Hall, Ohio
University, Athens, Ohio 45701, USA}

\email{uspensk@math.ohiou.edu}

\thanks{{\it 2000 Mathematics Subject Classification:}
Primary: 54D30. Secondary: 54C15, 22C05.}

\date{1 November 2003}

\keywords{Valdivia compact space, open map, retract, indecomposable group}

\begin{abstract} 
A compact space $K$ is {\em Valdivia compact} if it can be embedded
in a Tikhonov cube $I^A$ in such a way that the intersection $K\cap\Sigma$
is dense in $K$, where $\Sigma$ is the sigma-product (= the set of points 
with countably many non-zero coordinates). We show that there exists a compact
connected Abelian group of weight $\o_1$ which is not Valdivia compact, and deduce that 
Valdivia compact spaces are not preserved by open maps.
\end{abstract}

\maketitle

\setcounter{tocdepth}{1}

\section{Introduction} 
Let $\{X_\a:\a\in A\}$ be a family of spaces, and let a point $x_\a^*$
be given in each $X_\a$. 
The {\em sigma-product} $\Si$ of the family $\{X_\a\}$ with the base point
$\{x_\a^*\}$
is the subset of the product $\prod X_\a$ consisting of all points
$\{x_\a\}$ such that the set $\{\a\in A:x_\a\ne x_\a^*\}$ is countable.
A compact space $K$ is {\em Valdivia compact} if it can be embedded in
a Tikhonov cube $I^A$ in such a way that $K\cap \Si$ is dense in $K$, 
where $\Si$ is the sigma-product of intervals with the zero base point.
The class of Valdivia compact spaces is a natural extension of the class
of Corson compact spaces, which are defined as compact subspaces of 
sigma-products of intervals.
A compact space is Corson compact if and 
only if it is Valdivia compact and countably tight. While
Corson compact spaces are preserved by maps (we use the word ``map" to mean
a continuous map), this is not true for Valdivia compact spaces (however, 
if a Valdivia compact space is mapped onto a countably tight space $X$, 
then $X$ is Corson compact \cite{Kalenda7}).

The class of Valdivia compact spaces was introduced by Argyros, Mercourakis and Negrepontis in \cite{AMN}, they showed among others that these spaces admit `sufficiently many' retractions, which gives a projectional resolution of the identity on their spaces of continuous functions, see \cite{Valdivia}.
The name {\em Valdivia compact} was introduced by Deville and Godefroy in \cite{DG}. 
Valdivia compacta have been extensively studied by Kalenda; we refer to his survey article \cite{Kalenda}. Kalenda proved in \cite{Kalenda5} that an open image of a Valdivia compact space is Valdivia provided it contains a dense
set of $G_\delta$ points. The general question whether Valdivia compact spaces 
are preserved by 
open maps had remained open, see \cite{Kalenda}.
In the present paper we answer this question in the negative.

To this end, we construct a compact connected Abelian group $G$ which is not
Valdivia compact. Every compact Abelian group is a homomorphic image of a product
of compact metrizable groups. Indeed, in virtue of the Pontryagin duality this 
assertion is equivalent to the following: 
every Abelian group embeds into a direct sum of countable groups.
To see that this is true, note that every Abelian group embeds into a
divisible group \cite[Thm. 24.1]{F1}, and every Abelian divisible group is a direct
sum of countable groups \cite[Thm. 23.1]{F1}. 
Noting that: (1) any
continuous onto homomorphism between compact groups is open; (2) compact metric
spaces are Valdivia; (3) Valdivia compact spaces are preserved by products
-- we conclude that open maps do not preserve Valdivia compactness.

Our compact group $G$ is the Pontryagin dual of an uncountable indecomposable
torsion-free Abelian group $A$. The above argument can be made easier in this
case: $A$ is a subgroup of the vector space $A\otimes\Q$ over the field $\Q$
of rationals, hence $G=A^*$ is a homomorphic image of a power of the compact
metrizable group $\Q^*$. 

The class of Valdivia compact spaces is contained in a wider class
that was denoted by $\sR$ in \cite{BKT}. Recall the definition
of this class. A map $f:X\to Y$ is {\em right-invertible} if there exists a map
$g:Y\to X$ such that $fg:Y\to Y$ is the identity. A map $f:X\to Y$ is 
right-invertible if and only if it is homeomorphic to a retraction.

The class $\sR$ is defined as the smallest class containing all 
compact metric spaces which is closed 
under inverse limits of continuous transfinite sequences whose bonding mappings are 
right-invertible (in that case the limit projections are 
right-invertible as well \cite[Proposition~4.6]{BKT}). An inverse sequence
$\{X_\a;p_\a^\beta:\a<\beta<\kappa\}$ is {\em continuous} if for every limit
ordinal $\delta<\kappa$ the space $X_\delta$ is naturally homeomorphic to
$\liminv\{X_\a;p_\a^\beta:\a<\beta<\delta\}$.
To see that Valdivia compact spaces are in $\sR$, 
note that every Valdivia compact space of uncountable weight
is the inverse limit of a continuous transfinite sequence of 
Valdivia compacta of smaller weight whose all bonding maps are retractions, see \cite{AMN} or \cite[Thm. 3.6.2]{Kalenda}.

The compact group that we construct is not in the class $\sR$. Thus our example
shows that the image of a product of compact metric spaces under an open map
need not be in $\sR$. Our example has the smallest possible weight, namely $\o_1$. 
It has been proved in \cite{KM} that a 0-dimensional open image of
a Valdivia compact space is Valdivia if its weight does not exceed $\o_1$. 

Let us also mention that every compact group is a Dugundji space \cite{U1, U2}
and every 0-dimensional Dugundji space is Valdivia \cite{KM}. It is unknown 
whether the class of Valdivia compacta is stable under retractions; in \cite{KM} 
an affirmative answer is given in the case where the retract has weight $\o_1$.

\section{Proof of the main theorem}

\begin{thm}[Main Theorem]
There exists a compact connected Abelian group of weight $\o_1$ which is not in the class $\sR$
and hence is not Valdivia compact.
\label{main}
\end{thm}

We need some prerequisites.

\begin{prop} Let a compact space $X$ be in the class $\sR$. For every countable
Abelian group $G$ and every integer $n\ge0$ the cohomology group $H^n(X,G)$
is covered by its countable direct summands.
\label{p1}
\end{prop}

The cohomology theory that we use
is the \v Cech theory or any of its equivalents (Alexander -- Spanier, sheaves,
etc.), {\em not} the singular theory. 

\begin{proof}
The cohomology functor turns inverse limits of compact spaces into direct
limits of Abelian groups and right-invertible maps into left-invertible
homomorphisms, which are injective homomorphisms onto a direct summand.
If $X$ is compact metric, the group $H^n(X,G)$ is countable. (To see this, 
consider first the case of compact polyhedra, and then represent $X$ as
the limit of an inverse sequence of polyhedra.) It follows that 
the class $\sR'$ of all compact spaces $X$ for which the proposition holds 
contains compact metric spaces and is stable under
limits of continuous inverse sequences with right-invertible bonding maps
and projections. Since $\sR$ is the smallest class with these properties,
we have $\sR\sbs \sR'$.
\end{proof}

Let $\T=\{z\in \C:|z|=1\}$ be the circle group. For an Abelian group $A$
we denote by $A^*$ its  Pontryagin dual (= the group of all characters
$\chi: A\to \T$, considered as a compact group, see e.g. \cite[Ch. 6]{HR}).

The following proposition is well-known.

\begin{prop} For every torsion-free Abelian group $A$
there exists a natural isomorphism $\phi:A\to H^1(X, \Z)$, where
$X=A^*$.
\label{p2}
\end{prop}

\begin{proof} The homomorphism $\phi$
can be described as follows. Every $a\in A$
can be identified with a character $\chi_a:X\to\T$. Pick a generator
$u\in H^1(\T,\Z)$, and put $\phi(a)=H^1(\chi_a)(u)\in H^1(X,\Z)$. 

If $A=\Z^n$ and $X=A^*=\T^n$, 
it is clear that $\phi$ is an isomorphism. 
The general case follows by passing to limits:
every torsion-free Abelian group is the direct limit of finitely-generated
free groups; the Pontryagin duality turns direct limits of discrete groups
into inverse limits of compact groups; and the cohomology functor turns
inverse limits of compact spaces back to direct limits.
\end{proof}

{\it Proof of Theorem~\ref{main}.}
There exists a torsion-free Abelian group $A$ of cardinality $\o_1$ which is 
indecomposable \cite[Sections~88 and~89]{F},
that is, $A$ has no proper direct summands. Let $X=A^*$. We claim that
the compact group $X$ has the required properties.

The duals of torsion-free
discrete groups are connected \cite[Theorem 24.25]{HR}, so $X$ is connected.
According to Proposition~\ref{p2}, the cohomology group $H^1(X,\Z)$ is isomorphic
to $A$ and therefore indecomposable. Proposition~\ref{p1} implies that 
$X$ is not in the class $\sR$. \qed

\begin{corol}
There exists a compact metric space $K$ and an open onto map $f:K^{\o_1}\to X$
such that $X$ is not in the class $\sR$ and hence not Valdivia.
\end{corol}

\begin{proof}
We explained the construction in Section~1: Let $A$ and $X=A^*$ be as in the 
preceding proof. Embed $A$ into $A\otimes\Q=\Q^{(\o_1)}$. Passing to the duals,
we get a homomorphism of compact groups $K^{\o_1}\to X$, where $K=\Q^*$ 
(we consider $\Q$ as a discrete group). 
\end{proof}

\end{document}